\documentclass[11pt]{amsart}

\usepackage {amsfonts,amsthm}
\usepackage[english]{babel}

\theoremstyle{plain}
\newtheorem{Th}{Theorem}

\def\R{{\mathbb{R}}}
\def\Z{{\mathbb{Z}}}
\def\N{{\mathbb{N}}}
\def\s{\sigma}
\def\d{\delta}
\def\f{\varphi}

\def\e{{\varepsilon}}
\def\Lin{{\rm Lin}}

\begin{document}
\title
{Any discrete almost periodic set of finite type\\ is an ideal crystal}

\author{Favorov S.}

\address{Mathematical School, Kharkov National University, Swobody sq.4,
Kharkov, 61077 Ukraine}

 \email{Sergey.Ju.Favorov@univer.kharkov.ua}

\date{}

\begin{abstract}
 A discrete set in the Euclidian space is  almost periodic, if  the
 measure with the unite masses at points of the set is almost periodic in the weak sense.
We prove the following result: if $A$ is a discrete almost periodic set and the set
$A-A$ is discrete, then $A$ consists of a finite number of translates of a full rank
lattice.
\end{abstract}

\keywords{ almost periodic set, discrete set, ideal crystal, quasicrystals}

\subjclass{Primary: 52C23; Secondary: 11K70, 52C07}

\maketitle

The notion of a discrete almost periodic set in the complex plane is well known in the
theory of almost periodic holomorphic and meromorphic functions
(cf.~\cite{L},\cite{T}, \cite{FRR}, \cite{F}). Discrete almost periodic sets in the
$p$-dimensional Euclidian space appear later in the mathematical theory of
quasicrystals (cf.~\cite{La},\cite{M}). In this connection, in \cite{La} the question
(Problem 4.4) was raised, whether any discrete almost periodic set is just a finite
union of translations  of a full rank lattice in $\R^p$.  In \cite{FK} and \cite{FK1}
we showed that any almost periodic perturbation of a full rank lattice in $\R^p$ is an
almost periodic set. Hence, there exists a wide class of such sets. In the present
article we prove that discrete almost periodic sets with a certain additional property
have the form $L+F$ with a full rank lattice $L\subset\R^p$ and a finite set $F$.
\medskip

Let us recall some known definitions (see, for example, \cite{C}, \cite{R}).

\medskip
A continuous function $f(x)$ in $\R^p$ is {\it almost periodic}, if for any $\e>0$ the
set of $\e$-almost periods of $f$
 $$
  \{\tau\in\R^p:\,\sup_{x\in\R^p}|f(x+\tau)-f(x)|<\e\}
  $$
  is a relatively dense set in $\R^p$. The latter means that there is
  $R=R(\e)<\infty$ such that any ball of radius $R$ contains an $\e$-almost period of
  $f$.

A Borel measure $\mu$ in $\R^p$ is {\it almost periodic} if it is almost periodic in
the weak sense, i.e., for any continuous function $\f$ in $\R^p$ with a compact
support the convolution $\int\f(x+t)~d\mu(t)$
  is an almost periodic function in $x\in\R^p$.

A  discrete  set $A\subset\R^p$ is {\it almost periodic}, if its counting measure
 $\mu_A=\sum_{x\in A}\d_x$,
 where $\d_x$ is the unit mass at the point $x$, is almost periodic.

\medskip
 There is a geometric criterium for discrete sets to be
almost periodic.
\begin{Th}[\cite{FK}, Theorem 11]\label{FK}
 Suppose $(a_n)_{n\in\N}$ is an indexing of a discrete set $A\subset\R^p$. Then, the
 set $A$ is almost  periodic if and only if for each  $\e>0$ the set of $\e$-almost
 periods of $A$
 $$
\{\tau\in\R^p:\,\exists \hbox{\ a bijection}\quad \s:\,\N\to\N\quad\hbox{ such that
}\quad |a_n+\tau-a_{\sigma(n)}|<\e\quad \forall\,n\in\N\}
 $$
 is relatively dense in $\R^p$.
\end{Th}

It follows easily from this criterion  that the number of elements of $A$ in any ball
of radius $1$ is uniformly bounded.

\medskip
Following \cite{La}, we will say that a discrete set $A\subset\R^p$ is  of {\it finite
type}, if the set $A-A$ is discrete. A set $A\subset\R^p$ is an {\it ideal crystal},
if $A$ consists of a finite number of translates of a full rank lattice $L$, that is,
$A=L+F$, where $F$ is a finite set and $L$ is an additive discrete subgroup of $\R^p$
such that $\Lin_\R L=\R^p$.

In the present article we prove the following theorem.
\begin{Th}
  Any almost periodic discrete set of finite type is an ideal crystal.
\end{Th}

{\bf Proof}. From the definition of a relatively dense set it follows that there is
$D<\infty$ such that each ball with center at any point $a\in A$ and radius $D$
contains at least one point $b\in A,\,b\neq a$. Since the set $A-A$ is discrete, we
see that there is $\e>0$  such that $\e<\min\{1;|(a-b)-(c-d)|\}$ whenever
$a,\,b,\,c,\,d\in A$ and $|a-b|<D+1,\,|c-d|<D+1,\,a-b\neq c-d$.  In particular,
$\e<|a-b|$ whenever $a,\,b\in A$ and $a\neq b$.

Let $\tau\in\R^p$ be an arbitrary $(\e/2)$-almost period of $A$. Taking into account
our choice of $\e$, we see that there is a unique $c\in A$ that satisfies the
inequality $|a+\tau-c|<\e/2$. Clearly, $T=c-a$ is an $\e$-almost period of $A$. Let us
show that $T$ is actually a period of $A$.

Suppose that $b\in A$ such that $b\neq a$ and $|a-b|<D$. Since $T$ is an $\e$-almost
period of $A$, there exists a point $d\in A$ such that $|b+T-d|=|(a-b)-(c-d)|<\e$.
Since $|c-d|\le|a-b|+|b+T-d|<D+1$, we obtain $a-b=c-d$ and $d=b+T$. Repeating these
arguments for all $b\in A$ such that $|b-a|<D+1$, then, for all $b'\in A$ such that
$|b'-b|<D+1$. After a countable number of steps we obtain that $a+T\in A$ for all
$a\in A$.

For $p=1$, the conclusion of the theorem is evident. If $p>1$, for every
$j=1,\dots,p,$ take an $(\e/2)$-almost periods $\tau_j$ from the set
 $\{x\in\R^p:\,3p^2<|x|<(1+(2p)^{-2})~|\langle x,e_j\rangle|\}$,
where $e_j$ is the basis vector in $\R^p$. Since the corresponding period $T_j$
satisfies the inequality $|T_j-\tau_j|<1/2$, we get
 $$
 |T_j|<(1+2^{-1}p^{-2})~|\langle
T_j,e_j\rangle|\quad\hbox{and}\quad \max_{k\neq j}|\langle
T_j,e_k\rangle|<(p-1)^{-1}|\langle T_j,e_j\rangle|,\quad j=1,\dots,p.
 $$
 The latter inequalities imply that the determinant of the matrix
 $(\langle T_j,e_k\rangle)_{j,k=1}^p$ does not vanish, hence the vectors
$T_1,\dots,T_p$ are linearly independent. Consequently, the set $L=\{n_1T_1+\dots
+n_pT_p:\,n_1,\dots,n_p\in\Z\}$ is a full rank lattice. Next, the set $F=\{a\in A:
|a|<|T_1|+\dots+|T_p|\}$ is finite. All vectors $t\in L$ are periods of $A$, hence,
$L+F\subset A$. On the other hand, for each $a\in A$ there is $t\in L$ such that
$|a-t|<|T_1|+\dots+|T_p|$. The theorem is proved.

\end{document}